\documentclass[twocolumn]{autart}    

\def\baselinestretch{1.2}
\usepackage{tikz}

\usepackage{amsmath}
\usepackage{latexsym, amssymb}
\usepackage{graphicx}
\usepackage{amsmath, amsbsy}
\usepackage{amsopn, amstext}
\usepackage{hyperref}
\usepackage{cancel, color}
\usepackage{epstopdf}

\DeclareMathOperator{\rank}{rank}

\DeclareMathOperator{\Span}{Span}
\DeclareMathOperator{\Col}{Col}

\DeclareMathOperator{\lcm}{lcm}
\def\cal{\mathcal}

\def\diag{diag}
\def\ra{\rightarrow}

\def\d{\delta}

\def\D{\Delta}

\def\ot{\otimes}

\def\0{{\bf 0}}

\newcommand{\R}{{\mathbb R}}

\def\dsum{\mathop{\sum}\limits}
\def\thefootnote{*}

\newtheorem{dfn}[thm]{Definition}
\newtheorem{prp}[thm]{Proposition}
\newtheorem{exa}[thm]{Example}

\begin{document}

\begin{frontmatter}

\title{On Coset Weighted Potential Game\thanksref{footnoteinfo}}

\thanks[footnoteinfo]{This work was supported partly by National Natural Science Foundation (NNSF) of China under Grants 61773371. Corresponding author: Yuanhua Wang.}

\author{Yuanhua Wang}\dag \ead{wyh\_1005@163.com},~
\author{Daizhan Cheng}\ddag\ead{dcheng@iss.ac.cn},~

\address{\dag Business School, Shandong Normal University, Jinan250014, China\\
\ddag Key Laboratory of Systems and Control, Institute of Systems Science,
Chinese Academy of Sciences, Beijing 100190, China}

%
%
%
%

\begin{keyword}
coset weighted potential game, potential function, semi-tensor product of matrices, orthogonal decomposition.
\end{keyword}

\begin{abstract}
In this paper we first define a new kind of potential games, called  coset weighted potential game, which is a generalized form of weighted potential game. Using semi-tensor product of matrices, an algebraic method is provided to verify whether a finite game is a coset weighted potential game, and a simple formula is obtained to calculate the corresponding potential function. Then some properties of coset weighted potential games are revealed. Finally, by resorting to the vector space structure of finite games, a new orthogonal decomposition based on coset weights is proposed, the corresponding geometric and algebraic expressions of all the subspaces are given by providing their bases.
\end{abstract}

\end{frontmatter}

\section{Preliminaries}

 A finite normal game can be described by $G=(N,S,C)$, where $N=\{1,2,\cdots,n\}$ is the set of players; $S=\prod_{i=1}^nS_i$ is the strategy profile, and the set of strategies for player $i$ is $S_i=\{1,2,\cdots,k_i\}$. $S_{-i}=\prod_{j\neq i}S_j$ denotes the strategies of all players except the $i$-th one; $C=(c_1,\cdots,c_n)\in \R^n$ with $c_i:S\ra \R$ is the payoff function of player $i$. For statement ease, the set of finite games with $|N|=n$, $|S_i|=k_i$, $i=1,\cdots,n$, is denoted by ${\cal G}_{[n;k_1,\cdots,k_n]}$. As a special class of finite normal games, the potential game imposes restriction on the players' payoff functions. Potential game was first proposed by Rosenthal \cite{ros73}. Monderer and Shapley systematically investigated potential games and proved several useful properties in \cite{mon96b}, such as best response dynamics and fictitious play, converging to a Nash equilibrium, etc. Since then it has been applied to many engineering problems, including computer networks \cite{hei06}, distributed coverage of graphs \cite{zhu13}, and congestion control \cite{hao18}, etc. Several classes of potential games are described as follows.

 A function $P:S\ra \R$ is called an ordinal potential for $G$, if for any $x,y\in S_i$, and any $s_{-i}\in S_{-i}$, $i\in N$,
$$
c_i(x,s_{-i})-c_i(y,s_{-i})>0\Leftrightarrow P(x,s_{-i})-P(y,s_{-i})> 0,
$$
then $G$ is called an ordinal potential game.

In an ordinal potential game, only the signs of the difference in individual payoffs for each player, and the difference in  potential function, have to be the same. In fact, the really useful model in some physical applications is not the ordinal potential game, but the weighted (or exact) potential game \cite{mon96b}. Let $w=(w_i)_{i\in N}$ be a vector of positive weights, if there exists a function $P:S\ra \R$, called the weighted potential function, such that for any $x,y\in S_i$, and any $s_{-i}\in S_{-i}$, $i\in N$,
$$
c_i(x,s_{-i})-c_i(y,s_{-i})=w_i\left(P(x,s_{-i})-P(y,s_{-i})\right),
$$
then $G$ is called a weighted potential game. Especially, $G$ is called an exact potential game if $w_i=1$, $\forall i \in N$.


However, the weighted (or exact) potential games only cover a few class of games in practice. Moreover, a weighted potential game is essentially an exact potential game, because if we replace the payoffs $c_i$ by $c_i/w_i$, then a weighted potential game becomes an exact potential game. This fact stimulates us to find a more general weighted potential game, which is between the ordinal potential game and classical weighted potential game.

\begin{figure}[h!]
\begin{center}
\includegraphics[height=4.2cm,angle=0]{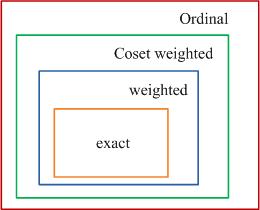}\label{fig1}\caption{Some classes of potential games}
\end{center}
\end{figure}

In this paper, we propose a new kind of weighted potential games, called the coset weighted potential game. Its relationship with classical kinds of potential games is depicted by Figure \ref{fig1}.
After a rigorous definition, we provided a simple method to verify whether a finite game is a coset weighted potential game. We show that though it is a generalization of classical weighted potential game and it can not be converted easily to exact potential game, it still has all the nice properties of classical (weighted) potential games. For instance, the existence of pure Nash equilibrium, the convergence to an equilibrium point under certain learning process, etc.

Another interesting topic for finite games is their vector space structure. In addition to (exact) potential games there are some other important kinds of finite games, which are necessary for investigating the vector space structure of ${\cal G}_{[n;k_1,\cdots,k_n]}$.

\begin{dfn}\label{d2.2.2} \cite{can11,liu15}  Let $G\in {\cal G}_{[n;k_1,\cdots,k_n]}$.
\begin{enumerate}
\item $G$ is called a non-strategic game if for any $x,y\in S_i$, and any $s_{-i}\in S_{-i}$,
$$
c_i(x,s_{-i})=c_i(y,s_{-i}),\quad i=1,\cdots,n.
$$
\item $G$ is called a harmonic game, if for any $s\in S$, and any $s_{-i}\in S_{-i}$,
$$
\dsum_{i=1}^n( c_i(s)-\frac{1}{k_i}\dsum_{x_i\in S_i}c_i(x_i,s_{-i}))=0.
$$
\item $G$ is called a pure harmonic game, if for any $s\in S$, and any $s_{-i}\in S_{-i}$,
$$
\begin{array}{llc}
\dsum_{i=1}^n c_i(s)=0;~~\dsum_{x\in S_i} c_i(x,s_{-i})=0,\quad i=1,\cdots,n.
\end{array}
$$
\end{enumerate}
\end{dfn}

Without the weights, by using the Helmholtz decomposition theorem, an orthogonal decomposition of ${\cal G}_{[n;k_1,\cdots,k_n]}$, briefly denoted by ${\cal G}$, was first proposed in \cite{can11}, which is described as follows.
\begin{align}\label{1.4}
{\cal G}=\rlap{$\underbrace{\phantom{\quad{\cal P}\quad\oplus\quad{\cal N}}}_{potential~ games}$}\quad{\cal P}\quad\oplus\quad
\overbrace{{\cal N}\quad\oplus\quad{\cal H}}^{harmonic~ games},
\end{align}
 where ${\cal P}$ is the subspace of pure potential games, ${\cal N}$ is the subspace of non-strategic games, and ${\cal H}$ is the subspace of pure harmonic games. An alternatively simplified approach was provided in \cite{che16b} to precisely express the bases of these orthogonal subspaces, by using the conventional inner product of Euclidean space. As a generalization of (\ref{1.4}), in this paper we are ready to proved a new orthogonal decomposition based on coset-depending weights for ${\cal G}$, which is described as follows.
\begin{align}\label{1.5}
{\cal G}=\rlap{$\underbrace{\phantom{\quad\quad{\cal P}^{cw}\quad\quad\quad\oplus{\cal N}}}_{coset~ weighted ~potential~ games~{\cal G}^{cw}_P}$}\quad{\cal P}^{cw}\quad\quad\quad\oplus
\overbrace{{\cal N}\quad\quad\oplus\quad\quad{\cal H}^{cw}}^{coset~weighted~ harmonic~games},
\end{align}
where ${\cal P}^{cw}$ is called the coset weighted pure potential subspace, and ${\cal H}^{cw}$ is called the coset weighted pure harmonic subspace. For each subspace, we give its geometric expression by providing the basis. Based on these bases, we also give an algebraic expression for each subspace, that is, the algebraic equation for the payoffs of the corresponding games to be satisfied. Meanwhile, some formulas are presented to calculate all the decomposed subspaces.

For statement ease, we first introduce some notations:
\begin{enumerate}

\item[$\bullet$]  $\scalebox{.9}{${\cal M}_{m\times n}$}$: the set of $m\times n$ real matrices.

\item[$\bullet$] $\scalebox{.9}{$\Col(M)$}$: the set of columns of $\scalebox{.9}{$M$}$. $\scalebox{.9}{$\Col_i(M)$}$: the $i$-th column of $\scalebox{.9}{$M$}$.

\item[$\bullet$] $\scalebox{.9}{${\cal D}_{k_i}:=\left\{1,2,\cdots,k_i\right\},\quad k_i\geq 2$}$.

\item[$\bullet$] $\d_n^i$: the $i$-th column of the identity matrix $\scalebox{.9}{$I_n$}$.

\item[$\bullet$] $\scalebox{.9}{$\D_n:=\{\d_n^i\vert i=1,\cdots,n\}$}$.

\item[$\bullet$] $\scalebox{.9}{${\bf 1}_{\ell}=(\underbrace{1,1,\cdots,1}_{\ell})^T$}$;
$\scalebox{.9}{${\bf 0}_{\ell}=(\underbrace{0,0,\cdots,0}_{\ell})^T$}$.

\item[$\bullet$] $\scalebox{.9}{${\bf 0}_{p\times q}$}$: a $p\times q$ matrix with zero entries.

\item[$\bullet$] $\scalebox{.9}{$[i,j]:=\{i,i+1,\cdots,j\}$}$, where $i,~j$ are integers and $i<j$.

\item[$\bullet$] A matrix $\scalebox{.9}{$L\in {\cal M}_{m\times n}$}$ is called a logical matrix
if the columns of $L$ are of the form
$\d_m^k$. Denote by $\scalebox{.9}{${\cal L}_{m\times n}$}$ the set of $m\times n$ logical
matrices. If $\scalebox{.9}{$L\in {\cal L}_{n\times r}$}$, by definition it can be expressed as
$\scalebox{.9}{$L=[\d_n^{i_1},\d_n^{i_2},\cdots,\d_n^{i_r}]$}$. It is briefly denoted as $\scalebox{.9}{$L=\d_n[i_1,i_2,\cdots,i_r]$}$.

\item[$\bullet$] $\scalebox{.9}{$\Span\{A_1,\cdots,A_s\}$}$: The subspace spanned by $\scalebox{.9}{$\{\Col(A_i)\;|\;i=1,\cdots,s\}$}$.

\item[$\bullet$] $\scalebox{.9}{$U\oplus V4$}$: orthogonal sum of two vector spaces, i.e., $u\perp v$, $\forall ~u\in U,~v\in V$.

\end{enumerate}

The semi-tensor product (STP) of matrices is a generalization of conventional matrix product, which is defined as follows \cite{che11}:

\begin{dfn}\label{d2.1.1}  Let $\scalebox{.9}{$M\in {\cal M}_{m\times n}$}$, $\scalebox{.9}{$N\in {\cal M}_{p\times q}$}$, and $t=\lcm\{n,p\}$ be the least common multiple of $n$ and $p$.
The STP of $\scalebox{.9}{$M$}$ and $\scalebox{.9}{$N$}$ is defined as
\begin{align}\label{2.1.1}
\scalebox{.9}{$M\ltimes N:= (M\otimes I_{t/n})(N\otimes I_{t/p})\in {\cal M}_{mt/n\times qt/p},$}
\end{align}
where $\otimes$ is the Kronecker product.
\end{dfn}
The STP keeps all the properties of the conventional matrix product. Hence we can omit the symbol $\ltimes$ mostly.
This method has been widely used to study the logical dynamic systems \cite{che09,gy16,lif2017,lih16,luj17,wu15}, and game theory \cite{che15,guo13}, etc. Next we give some properties of STP used in this paper.
\begin{prp}\label{p2.1} Let $\scalebox{.9}{$X\in \R^n$}$ be a column and $\scalebox{.9}{$M$}$ be a matrix. Then \scalebox{.9}{$X\ltimes M=\left(I_n\otimes M\right)X.$}
\end{prp}
\begin{prp}\label{p2.2} Let $\scalebox{.9}{$X\in \D_p$}$ and define a power reducing matrix \scalebox{.9}{$O^R_p:=\d_{p^2}[1,p+2,2p+3,\cdots,p^2]\in {\cal L}_{p^2\times p}$}. Then \scalebox{.9}{$X^2=O^R_pX.$}
\end{prp}

To use matrix expression for finite games, we identify each strategy $j\in {\cal D}_{k_i}$ by $\d_{k_i}^j$, that is, $j\sim \d_{k_i}^j$, $j=1,\cdots,k_i$, then $S_i\sim \D_{k_i}$, $i=1,\cdots,n$. It follows that the payoff functions can be expressed as
\begin{align}\label{2.2.2}
c_i(x_1,\cdots,x_n)=V^c_i\ltimes_{j=1}^nx_j,\quad i=1,\cdots,n,
\end{align}
where $V^c_i\in \R^{k}$ ($k=\prod_{i=1}^nk_i$) is a row vector, called the structure vector of $c_i$. Define the structure vector of a given game $G$ as
\begin{align}\label{2.2.3}
V_G=[V^c_1,V^c_2,\cdots,V^c_n]\in \R^{nk}.
\end{align}
It is clear that ${\cal G}$ has a natural vector space structure as ${\cal G}\sim \R^{nk}.$ For a given game $G\in {\cal G}$, its structure vector $V_G$ completely determines $G$. So the vector space structure is very natural and reasonable.

The rest of this paper is organized as follows: In Section 2 we give an algebraic verified method for coset weighted potential games. The results are used to verify the two-player Boolean game. Moreover, we present some important properties of coset weighted potential games. Section 3 derive a new orthogonal decomposition of finite games based on coset-depending weights. The geometric and algebraic expressions of all the subspaces are obtained by providing their bases. Based on these bases, some numerical formulas are provided for calculating all the decomposed components. Section 4 is a conclusion.

\section{Algebraic Verification of coset weighted potential games}
\subsection{Coset weighted potential equation}

We define coset weighted potential games as follows.
\begin{dfn}\label{d1.1}  A finite game $G=(N,S,C)$ is called a coset weighted potential game, if there exists a function $P:S\ra \R$, called the coset weighted potential function, and a set of weights $w_i(s_{-i})$ depending on $s_{-i}$, such that for any $x,y\in S_i$, and any $s_{-i}\in S_{-i}$, $i\in N$,
\begin{align}\label{1.3}
\scalebox{.9}{$c_i(x,s_{-i})-c_i(y,s_{-i})=w_i(s_{-i})\left(P(x,s_{-i})-P(y,s_{-i})\right).$}
\end{align}
\end{dfn}
Obviously, (\ref{1.3}) is equivalent to that there exists a function $d_i$, which is independent of $x\in S_i$, such that for any $x\in S_i$ and any $s_{-i}\in S_{-i}$,
\begin{align}\label{3.2}
\scalebox{.9}{$c_i(x,s_{-i})-w_i(s_{-i})P(x,s_{-i})=d_i(s_{-i}).$}
\end{align}
Using (\ref{2.2.2}), we express (\ref{3.2}) in its vector form as
\begin{align}\label{3.3}
\scalebox{.9}{$V^c_i\ltimes_{j=1}^nx_j-V^w_i\ltimes_{j\neq i}x_jV^P\ltimes_{j=1}^nx_j=V^d_i\ltimes_{j\neq i}x_j,$}
\end{align}
where \scalebox{.9}{$V^c_i,~V^P\in \R^{k}$}, and \scalebox{.9}{$V^w_i\in \R_+^{k/{k_i}}$, $V^d_i\in \R^{k/{k_i}}$} are the row vectors. Now verifying whether $G$ is coset weighted potential is equivalent to checking whether the solution of (\ref{3.3}) for unknown vectors $V^P$ and $V^d_i$ exists. Define a matrix operator as
\begin{align}\label{3.3.1}
\scalebox{.9}{$E_i:=I_{k^{[1,i-1]}}\otimes {\bf 1}_{k_i}\otimes I_{k^{[i+1,n]}}\in {\cal M}_{k\times k/k_i},~ i=1,\cdots,n,$}
\end{align}
where
$$
\scalebox{.9}{$k^{[p,q]}:=
\begin{cases}
\prod_{j=p}^qk_j,\quad q\geq p,\\
1,\quad q<p,
\end{cases}$}
$$
then (\ref{3.3}) becomes
$$
\scalebox{.9}{$V^c_i\ltimes_{j=1}^nx_j-V^w_iE_i^T \ltimes_{j=1}^nx_j V^P\ltimes_{j=1}^nx_j=V^d_iE_i^T \ltimes_{j=1}^nx_j.$}
$$
Using Proposition \ref{p2.1} and \ref{p2.2}, we have
$$
\scalebox{.9}{$V^c_i\ltimes_{j=1}^nx_j-V^w_iE_i^T ( I_k\otimes V^P)O_k^R\ltimes_{j=1}^nx_j=V^d_iE_i^T \ltimes_{j=1}^nx_j.$}
$$
It follows that
\begin{align}\label{3.4}
\scalebox{.9}{$V^c_i-V^w_iE_i^T (I_k\otimes V^P)O_k^R=V^d_iE_i^T,~i=1,\cdots,n.$}
\end{align}
Next we give a simple lemma.
\begin{lem}\label{l3.1} Let $X,~Y\in \R^n$ be two rows, then
$$
\scalebox{.9}{$X(I_n\otimes Y)=Y(X\otimes I_n).$}
$$
\end{lem}
 \noindent{\bf Proof.} Set $\scalebox{.9}{$X=[x_1,\cdots,x_n]$}$ and $\scalebox{.9}{$Y=[y_1,\cdots,y_n]$}$, a straight forward calculation shows that
 $$
 \begin{array}{llc}
 \scalebox{.9}{$X(I_n\otimes Y)=[x_1,x_2,\cdots,x_n](I_n\otimes [y_1,y_2,\cdots,y_n])$}\\
                \scalebox{.9}{$=[x_1y_1,x_1y_2,\cdots,x_1y_n,\cdots,x_ny_1,\cdots,x_ny_n],$}
 \end{array}
 $$
and
$$
\begin{array}{llc}
\scalebox{.9}{$ Y(X\otimes I_n)=[y_1,y_2,\cdots,y_n]([x_1,x_2,\cdots,x_n]\otimes I_n)$}\\
                \scalebox{.9}{$=[x_1y_1,x_1y_2,\cdots,x_1y_n,\cdots,x_ny_1,\cdots,x_ny_n].$}
\end{array}
$$
Hence, $\scalebox{.9}{$X(I_n\otimes Y)=Y(X\otimes I_n).$}$
\hfill $\Box$

Using Lemma \ref{l3.1}, (\ref{3.4}) becomes
\begin{align}\label{3.5}
\scalebox{.9}{$V^P (V^w_iE_i^T \otimes I_k)O_k^R=V^c_i-V^d_iE_i^T,~i=1,\cdots,n.$}
\end{align}
Since \scalebox{.9}{$V^w_i\in \R_+^{k/{k_i}}$}, then \scalebox{.9}{$V^w_iE_i^T \in \R_+^k$}. Denote $V^w_iE_i^T=[w_i^1,w_i^2,\cdots,w_i^k]$, according to Definition \ref{d2.1.1}, we have
$$
\scalebox{.9}{$(V^w_iE_i^T \otimes I_k)O_k^R=V^w_iE_i^T \ltimes O_k^R=\diag(w_i^1,w_i^2,\cdots,w_i^k).$}
$$
Denote \scalebox{.9}{$\Lambda_i=V^w_iE_i^T \ltimes O_k^R$, $i=1,2,\cdots,n$}. Obviously, the diagonal matrix $\Lambda_i$ is reversible. Solving $V^P$ from the first equation of (\ref{3.5}) yields
$$
\scalebox{.9}{$V^P$}\scalebox{.9}{$=(V^c_1-V^d_1E_1^T)(V^w_1E_1^T \ltimes O_k^R)^{-1}=(V^c_1-V^d_1E_1^T)\Lambda_1^{-1}.$}
$$
Plugging it into the rest equations of (\ref{3.5}) yields
$$
\scalebox{.9}{$(V^c_1-V^d_1E_1^T)\Lambda_1^{-1}\Lambda_i=V^c_i-V^d_iE_i^T,~i=2,\cdots,n.$}
$$
It follows that
$$
\scalebox{.9}{$(V^c_1-V^d_1E_1^T)\Lambda_1^{-1}=(V^c_i-V^d_iE_i^T)\Lambda_i^{-1},~ i=2,\cdots,n.$}
$$
Taking transpose, we have
$$
\scalebox{.9}{$\Lambda_1^{-1}\left[(V^c_1)^T-E_1(V^d_1)^T\right]=\Lambda_i^{-1}\left[(V^c_i)^T-E_i(V^d_i)^T\right].$}
$$
It can be rewritten as
$$
\scalebox{.9}{$-\Lambda_1^{-1}E_1(V^d_1)^T+\Lambda_i^{-1}E_i(V^d_i)^T=\Lambda_i^{-1}(V^c_i)^T-\Lambda_1^{-1}(V^c_1)^T.$}
$$
Since $\Lambda_i$ are all diagonal matrices, they are mutually commutative. For the above equation, we first left multiply both sides by $\Lambda_1$ and $\Lambda_i$, we have
\begin{align}\label{3.7}
\begin{array}{llc}
\scalebox{.9}{$-\Lambda_iE_1(V^d_1)^T+\Lambda_1E_i(V^d_i)^T=\Lambda_1(V^c_i)^T-\Lambda_i(V^c_1)^T,$}\\
~~~~~~~~~~~~~~~~~~~~~~~~~~~~~~~~~\quad \quad i=2,\cdots,n.
\end{array}
\end{align}
Define $\xi^w_i:=\left(V^d_i\right)^T\in \R^{k/{k_i}},~i=1,\cdots,n,$ and
$$
\scalebox{.9}{$b^w_i:=\Lambda_1(V^c_i)^T-\Lambda_i(V^c_1)^T\in \R^{k},~ i=2,\cdots,n.$}
$$
(\ref{3.7}) can be expressed as a linear system:
\begin{align}\label{3.8}
\Psi_w\xi^w=b^w,
\end{align}
where $\xi^w=[\xi_1^w,\xi_2^w,\cdots,\xi_n^w]^T$, $b^w=[b_2^w,b_3^w,\cdots,b_n^w]^T$, and
$$
\Psi_w=\begin{bmatrix}
\begin{smallmatrix}
-\Lambda_2E_1&\Lambda_1E_2&0&\cdots&0\\
-\Lambda_3E_1&0&\Lambda_1E_3&\cdots&0\\
\vdots&\vdots&\vdots&\ddots&\vdots\\
-\Lambda_nE_1&0&0&\cdots&\Lambda_1E_n
\end{smallmatrix}
\end{bmatrix}.
$$
Eq.(\ref{3.8}) is called the coset weighted potential equation and $\Psi_w$ is called the coset weighted potential
matrix. Then we have the following result.
\begin{thm}\label{t3.1} A finite normal game $G\in {\cal G}$ is a coset weighted potential game with a set of coset-depending weights $w_i(s_{-i})>0$, if and only if Eq.(\ref{3.8}) has solutions. Moreover, the coset weighted potential is
\begin{align}\label{3.8.0}
\scalebox{.9}{$V^P=(V^c_1-V^d_1E_1^T)\Lambda_1^{-1}.$}
\end{align}
\end{thm}

\begin{rem}\label{r3.0.1} In \cite{che14c}, the potential matrix $\Psi$ depends on $n$ and $k_i=|S_i|$, while $b$ depends on the payoffs, so only the payoffs determines whether a game is an exact potential game. However, the matrix $\Psi_w$ in (\ref{3.8}) not only depends on $n$ and $k_i=|S_i|$, but coset-depending weights $w_i(s_{-i})$, while $b^w$ depends on the payoffs and $w_i(s_{-i})$, then if $G$ is not an exact potential game, choosing its coset-depending weights can make it a coset weighted potential game.
\end{rem}

\subsection{Two-player Boolean game}
As a simple application, we consider a game of two players with two strategies for each player, which is called a two-player Boolean game \cite{chengd18}. Denote ${\cal G}_{[2;2,2]}$ as the set of two-player Boolean games.
\begin{exa}\label{ex1} Consider a two-player Boolean game $G\in {\cal G}_{[2;2,2]}$. Its payoffs can be expressed in Table \ref{Tab1}.
\begin{table}[!htbp] 
\caption{Payoffs of a two-player Boolean game}
\label{Tab1}
\doublerulesep 0.1pt
\begin{tabular}{|c||c|c|}
\hline $P_1\backslash P_2$&$1$&$2$\\
\hline $1$&$(a,~e)$&$(b,~f)$\\
\hline $2$&$(c,~g)$&$(d,~h)$\\
\hline
\end{tabular}
\end{table}
Using the potential equation in \cite{wang17}, it is easy to verify that $G$ is a weighted potential game when its payoffs satisfy $(a-b-c+d)(e-f-g+h)>0$, otherwise it is not. If yes, its weights satisfy
$$
\frac{w_1}{w_2}=\frac{a-b-c+d}{e-f-g+h}.
$$
Particularly, $G$ becomes an exact potential game when $a-b-c+d=e-f-g+h$. Consider $a=-1$, $b=2$, $c=0$, $d=3$, $e=3$, $f=3$, $g=5$, $h=4$. it follows that $(a-b-c+d)=0$, $(e-f-g+h)\neq 0$. Obviously, it is not a weighted potential game. However, we can choose suitable coset-depending weights $w_i(s_{-i})$ for player $i$, $i=1,2,$ which makes $G$ a coset weighted potential game.

Assume $V^w_i=[\alpha_i,\beta_i]$, $\alpha_i,\beta_i>0$, $i=1,2$, then we have
$$
\begin{array}{llc}
\scalebox{.9}{$\Lambda_1=V^w_1E_1^T \ltimes O_k^R=\diag(\alpha_1,\beta_1,\alpha_1,\beta_1),$}\\
\scalebox{.9}{$\Lambda_2=V^w_2E_2^T \ltimes O_k^R=\diag(\alpha_2,\alpha_2,\beta_2,\beta_2).$}
\end{array}
$$
According to (\ref{3.8}), we obtain that
$$
\begin{bmatrix}
\begin{smallmatrix}
-\alpha_2&0&\alpha_1&0\\
0&-\alpha_2&\beta_1&0\\
-\beta_2&0&0&\alpha_1\\
0&-\beta_2&0&\beta_1
\end{smallmatrix}
\end{bmatrix}
\begin{bmatrix}\begin{smallmatrix}\xi_1^w\\\xi_2^w\end{smallmatrix}\end{bmatrix}=\begin{bmatrix}\begin{smallmatrix}3\alpha_1+\alpha_2\\3\beta_1-2\alpha_2 \\ 5\alpha_1 \\
4\beta_1-3\beta_2\end{smallmatrix}\end{bmatrix}.
$$
Choose $\alpha_1=1$, $\beta_1=2$, $\alpha_2=3$, $\beta_2=2$, the above equation has solutions and one of solutions can be solved out as
$$
\scalebox{.9}{$\xi_1^w=\left(V^d_i\right)^T=[-2.5,-1]^T.$}
$$
Using (\ref{3.8.0}), the coset weighted potential is calculated as
$$
\scalebox{.9}{$V^P=(V^c_1-V^d_1E_1^T)\Lambda_1^{-1}=[ 1.5 ,~1.5, ~2.5,~ 2].$}
$$
Hence, by choosing coset-depending weights $w_i(s_{-i})$, the two-player Boolean game $G$ becomes a coset weighted potential game.
\end{exa}

From Example \ref{ex1}, it is shown that the coset weighted potential games is more general than the weighted potential games. Moreover, using (\ref{3.8}), the coset weighted potential equation can be expressed as
\begin{align}\label{e3.9}
\begin{bmatrix}\begin{smallmatrix}
-\alpha_2&0&\alpha_1&0\\
0&-\alpha_2&\beta_1&0\\
-\beta_2&0&0&\alpha_1\\
0&-\beta_2&0&\beta_1
\end{smallmatrix}\end{bmatrix}\begin{bmatrix}\begin{smallmatrix}\xi_1^w\\\xi_2^w\end{smallmatrix}\end{bmatrix}
=\begin{bmatrix}\begin{smallmatrix}\alpha_1e-\alpha_2a\\ \beta_1f-\alpha_2 b \\ \alpha_1g- \beta_2c \\ \beta_1h-\beta_2 d\end{smallmatrix}\end{bmatrix}.
\end{align}
Since $\rank(\Psi_w)=\rank(\Psi_w,b^w)$, by the straightforward computation, we have the following result.

\begin{prp}\label{p3.1} A two-player Boolean game $G\in {\cal G}_{[2;2,2]}$, with coset-depending weights $w_i(s_{-i})=[\alpha_i,\beta_i]\ltimes_{j\neq i}x_j$, $\alpha_i,\beta_i>0$, $i=1,2$, is a coset weighted potential game,
if and only if Eq. (\ref{e3.9}) has solutions, that is, the payoffs and coset-depending weights satisfy
$$
\scalebox{.9}{$\frac{1}{\alpha_1}(c-a)+\frac{1}{\alpha_2}(e-f)+\frac{1}{\beta_1}(b-d)+\frac{1}{\beta_2}(h-g)=0.$}
$$
Moreover, assume $[A, B, C,D]^T$ is a particular solution of (\ref{e3.9}), the coset weighted potential function can be obtained as
\begin{align}\label{3.10}
\scalebox{.9}{$P(x_1,\cdots,x_n)=V^P\ltimes_{j=1}^nx_j+c_0,\quad \forall c_0\in \mathbb{R},$}
\end{align}
where
$$
\begin{array}{llc}
V^P&=\left([a,b,c,d]-[A,B,A,B]\right)\Lambda_1^{-1}\\
&=\left[\frac{a-A}{\alpha_1},\frac{b-B}{\beta_1},\frac{c-A}{\alpha_1},\frac{d-B}{\beta_1}\right].
\end{array}
$$
\end{prp}

\subsection{Properties of coset weighted potential games}
For an exact potential game $G$, it was proved in \cite{mon96b} that the potential function $P$ is unique up to a constant number. That is, if $P_1$ and $P_2$ are two potential functions, then $P_1-P_2=c_0 \in \R$. A coset weighted potential game has the same property, which is described as follows.
\begin{prp}\label{p3.2} Consider a coset weighted potential game $G$. Let $P_1$ and $P_2$ are two coset weighted potential functions for $G$, then there exists a constant $c$ such that for every $s\in S$,
\begin{align}\label{e3.11}
\scalebox{.9}{$P_1(s)-P_2(s)=c\in \R.$}
\end{align}
\end{prp}
\noindent{\bf Proof.} From (\ref{1.3}) and (\ref{3.2}), if $P_1$ and $P_2$ are two potential functions for a coset weighted potential game $G$, then we have
$$
\begin{array}{llc}
\scalebox{.9}{$c_i(s)-w_i(s_{-i})P_1(s)=d_i(s_{-i}),$}\\
\scalebox{.9}{$c_i(s)-w_i(s_{-i})P_2(s)=d'_i(s_{-i}),$}
\end{array}
$$
Set \scalebox{.9}{$c=P_1(s)-P_2(s)$,} then
$$
c=\frac{d'_i(s_{-i})-d_i(s_{-i})}{w_i(s_{-i})}.
$$
$d'_i(s_{-i})$, $d_i(s_{-i})$ and $w_i(s_{-i})$ are all independent of $x\in S_i$, so $c$ is independent of $x\in S_i$. But player $i$ is arbitrary, hence, $c$ is a constant.
\hfill $\Box$

 According to Definition \ref{d1.1}, for a fixed coset-depending weights $w_i(s_{-i})$, it is easy to see that, in a coset weighted potential game, any strategy profile $s\in S$ maximizing the potential function $P$ is a pure strategy equilibrium. Hence, we have the following property.
\begin{prp}\label{p3.3} Consider a coset weighted potential game $G$ with fixed coset-depending weights $w_i(s_{-i})>0$. The game $G$ possesses at least one pure Nash equilibrium.
\end{prp}
Because of the existence of pure Nash equilibrium, there are many learning algorithms which lead a coset weighted potential game to a pure Nash equilibrium. For instance, it is easily proved that the Myopic Best Response Adjustment \cite{you93}, Fictitious Play \cite{mon96a}, etc, will all guarantee the convergence of a coset weighted potential game to one of pure Nash equilibria.

\section{Decomposition of finite games with coset-depending weights}

In this section, we respectively discuss the geometric and algebraic expressions of all the subspaces in (\ref{1.5}) by providing their bases. Based on these bases, (\ref{1.5}) is proved to be hold and some formulas are provided for calculating all the decomposed components.
\subsection{Subspace of coset weighted potential games ${\cal G}^{cw}_P$}

According to Theorem \ref{t3.1}, we can derive that $G\in {\cal G}$ is a coset weighted potential game with a set of coset-depending weights $w_i(s_{-i})$, if and only if
\begin{align}\label{3.10}
\scalebox{.9}{$b^w\in \Span(\Psi_w).$}
\end{align}
Observing that in (\ref{3.10}) we have freedom to choose arbitrarily $V^c_1$, then (\ref{3.10}) can be rewritten as
$$
\begin{bmatrix}\begin{smallmatrix}
(V^c_1)^T\\
\Lambda_1(V^c_2)^T-\Lambda_2(V^c_1)^T\\
\vdots\\
\Lambda_1(V^c_n)^T-\Lambda_n(V^c_1)^T
\end{smallmatrix}\end{bmatrix}\in \scalebox{.9}{$\Span(E_{cw}^e),$}
$$
where $\scalebox{.9}{$E_{cw}^e$}=\begin{bmatrix}\begin{smallmatrix} I_k&0\\0&\Psi_w\end{smallmatrix}\end{bmatrix}.$ It is equivalent to
$$
\begin{bmatrix}\begin{smallmatrix}
I_k&0&\cdots&0\\
-\Lambda_2&\Lambda_1&\cdots&0\\
\vdots&\vdots&\ddots&\vdots\\
-\Lambda_n&0&\cdots&\Lambda_1\\
\end{smallmatrix}\end{bmatrix}
\begin{bmatrix}\begin{smallmatrix}
(V^c_1)^T\\
(V^c_2)^T\\
\vdots\\
(V^c_n)^T\\
\end{smallmatrix}\end{bmatrix}\in \scalebox{.9}{$\Span(E_{cw}^e).$}
$$
It follows that $\scalebox{.9}{$V_G^T\in \Span(E_{cw}^P)$}$, where
$$
\begin{array}{lll}
\scalebox{.9}{$E_{cw}^P$}&=\begin{bmatrix}\begin{smallmatrix}
I_k&0&\cdots&0\\
-\Lambda_2&\Lambda_1&\cdots&0\\
\vdots&\vdots&\ddots&\vdots\\
-\Lambda_n&0&\cdots&\Lambda_1\\
\end{smallmatrix}\end{bmatrix}^{-1}\scalebox{.9}{$E_{cw}^e$}\\
~\\
&=\begin{bmatrix}\begin{smallmatrix}
I_k&0&0&\cdots&0\\
\Lambda_2\Lambda_1^{-1}&-\Lambda_1^{-1}\Lambda_2E_1&E_2&\cdots&0\\
\vdots&\vdots&~&\ddots&\vdots\\
\Lambda_n\Lambda_1^{-1}&-\Lambda_1^{-1}\Lambda_nE_1&0&\cdots&E_n\end{smallmatrix}\end{bmatrix}
~~\\
~\\
&=\scalebox{.9}{$D_w$}\begin{bmatrix}\begin{smallmatrix}
\Lambda_1&0&0&\cdots&0\\
\Lambda_2&-\Lambda_2E_1&\Lambda_1E_2&\cdots&0\\
\vdots&\vdots&~&\ddots&\vdots\\
\Lambda_n&-\Lambda_nE_1&0&\cdots&\Lambda_1E_n\end{smallmatrix}\end{bmatrix}
\end{array}
$$
where $\scalebox{.9}{$D_w=\diag(\Lambda_1^{-1},\Lambda_1^{-1},\cdots,\Lambda_1^{-1}).$}$ Assume $V_i^w={\bf 1}^T_{k/{k_i}}$ for any $i$, then the coset weighted potential games become the (exact) potential games. Similar to the arguments in \cite{che14c} and \cite{wang17}, we construct $E_{cw}^{P^0}$ from $E_{cw}^P$ via deleting the last column of $\Lambda_1E_n$, then $E_{cw}^{P^0}$ has full column rank. Hence, we have the following result.
\begin{thm}\label{t3.2} The subspace of coset weighted potential games is
\begin{align}\label{3.11}
\scalebox{.9}{${\cal G}_{P}^{cw}=\Span(E_{cw}^P),$}
\end{align}
which has $\scalebox{.9}{$\Col(E_{cw}^{P^0})$}$ as its basis.
\end{thm}
\begin{rem}\label{r3.31} The equation (\ref{1.3}) provides the algebraic condition for the payoff functions to satisfy, so (\ref{1.3}) is called the algebraic expression of coset weighted potential games. Moreover, (\ref{3.11}) is called the geometric expression of coset weighted potential games, because it gives the basis of the corresponding subspace.
\end{rem}

\subsection{Subspace of coset weighted pure potential games ${\cal P}^{cw}$}
Define
\begin{align}\label{3.14}
\scalebox{.9}{$\widetilde{E}_{cw}^P:$}=
\begin{bmatrix}\begin{smallmatrix}
\Lambda_1&\Lambda_1E_1&0&\cdots&0\\
\Lambda_2&0&\Lambda_1E_2&\cdots&0\\
\vdots&\vdots&\vdots&\ddots&\vdots\\
\Lambda_n&0&0&\cdots&\Lambda_1E_n\end{smallmatrix}\end{bmatrix}.
\end{align}
Compared (\ref{3.14}) and $E_{cw}^P$, we can verify that
$$
\scalebox{.9}{${\cal G}_{P}^{cw}=\Span(E_{cw}^P)=\Span(\widetilde{E}_{cw}^P).$}
$$

The subspace of non-strategic games \cite{che16b} is defined as ${\cal N}:=\Span(B^N),$ where
\begin{align}\label{3.15}
\scalebox{.9}{$B^N$}=\begin{bmatrix}\begin{smallmatrix}
E_1&0&\cdots&0\\
\vdots&~&\ddots&\vdots\\
0&0&\cdots&E_n
\end{smallmatrix}\end{bmatrix}.
\end{align}

Similar to the argument in \cite{che16b}, we define
\begin{align}\label{3.16}
\scalebox{.9}{$B_{cw}^P$}=\begin{bmatrix}\begin{smallmatrix}
\Lambda_1-\frac{1}{k_1}\Lambda_1E_1E_1^T\\
\Lambda_2-\frac{1}{k_2}\Lambda_2E_2E_2^T\\
\vdots\\
\Lambda_n-\frac{1}{k_n}\Lambda_nE_nE_n^T
\end{smallmatrix}\end{bmatrix}\in\scalebox{.9}{$ {\cal M}_{nk\times k}.$}
\end{align}
According to (\ref{3.15}) and (\ref{3.16}), it is easy to verify that $\scalebox{.9}{${\cal G}^{cw}_P=\Span\left\{B_{cw}^P,B^N\right\}$}$. Moreover, we can verify that $\scalebox{.9}{$\left(B_{cw}^P\right)^TB^N=0$}$. Hence, we have an orthogonal decomposition as $\scalebox{.9}{${\cal G}^{cw}_P=\Span\left\{B_{cw}^P\right\} \oplus {\cal N}.$}$
Obviously, the coset weighted pure potential subspace can be expressed as
\begin{align}\label{e3.16}
\scalebox{.9}{${\cal P}^{cw}:=\Span\left\{B_{cw}^P\right\}.$}
\end{align}
Since $\scalebox{.9}{$\dim({\cal P}^{cw})=k-1$}$, and $\scalebox{.9}{$B_{cw}^P{\bf 1}_k={\bf 0}_{nk}$}$, similar to the argument in \cite{che16b}, we can delete any one column of $B_{cw}^P$, say, the last column, and denote the remaining matrix by $\scalebox{.9}{$B_{cw}^{P^0}$}$, then we have
\begin{align}\label{e3.61}
\scalebox{.9}{${\cal P}^{cw}:=\Span\left\{B_{cw}^P\right\}=\Span\left\{B_{cw}^{P^0}\right\},$}
\end{align}
where $\scalebox{.9}{$\Col(B_{cw}^{P^0})$}$ is a basis of $\scalebox{.9}{${\cal P}^{cw}$}$.

According to (\ref{3.16}), we have the following result.
\begin{thm}\label{t3.4} Consider $G\in {\cal G}$. The following three statements are equivalent.
\begin{enumerate}
\item $G$ is a coset weighted pure potential game.
\item there exists a function $P:S\ra \R$ and a set of coset-depending weights $w_i(s_{-i})>0$, such that for any $s_{-i}\in S_{-i}$,
\begin{align}\label{e3.31}
\scalebox{.9}{$c_i(s)=w_i(s_{-i})P(s)-w_i(s_{-i})\sum_{x\in S_i}P(x,s_{-i}).$}
\end{align}
\item there exists a function $P:S\ra \R$ and a set of coset-depending weights $w_i(s_{-i})>0$, such that for any $s_{-i}\in S_{-i}$, and $x,~y\in S_i$,
\begin{align}
\label{e3.32}&\scalebox{.9}{$\sum_{x\in S_i}c_i(x,s_{-i})=0,\quad \forall s_{-i}\in S_{-i};$}\\
\label{e3.33}&\scalebox{.9}{$c_i(x,s_{-i})-c_i(y,s_{-i})=w_i(s_{-i})\left(P(x,s_{-i})-P(y,s_{-i})\right).$}
\end{align}
\end{enumerate}
\end{thm}
\noindent{\bf Proof.} $1\Rightarrow 2:$ According to (\ref{e3.16}), if $G$ is a coset weighted pure potential game, there exists a column $\gamma \in \R^k$, such that $\scalebox{.9}{$V_G^T=B_{cw}^P \gamma$}$. Set $X=\ltimes_{j=1}^nx_j$, define $\scalebox{.9}{$P(s)=\gamma^T\ltimes_{j=1}^nx_j=\gamma^TX$}$ , then we have
$$
\begin{array}{llc}
\scalebox{.9}{$c_i(s)=V^c_i\ltimes_{j=1}^nx_j=\gamma^T(\Lambda_i-\frac{1}{k_i}\Lambda_iE_iE_i^T)X$}\\
\scalebox{.9}{$=\gamma^T(V^w_iE_i^T \otimes I_k)O_k^RX-\frac{1}{k_i}\gamma^T(V^w_iE_i^T \otimes I_k)O_k^RE_iE_i^TX$}\\
\scalebox{.9}{$=V^w_iE_i^T (I_k\otimes \gamma^T)O_k^RX-\frac{1}{k_i}V^w_iE_i^T (I_k\otimes \gamma^T)O_k^RE_iE_i^TX$}\\
\scalebox{.9}{$=V^w_iE_i^TX\gamma^TX-\frac{1}{k_i}V^w_iE_i^TE_iE_i^TX\gamma^TE_iE_i^TX$}\\
\scalebox{.9}{$=w_i(s_{-i})P(s)-\frac{1}{k_i}V^w_iE_i^T\ltimes_{j=1}^{i-1}x_j\ltimes {\bf 1}_{k_i}\ltimes_{j=i+1}^nx_j$}\\
\scalebox{.9}{$\gamma^T\ltimes_{j=1}^{i-1}x_j\ltimes {\bf 1}_{k_i}\ltimes_{j=i+1}^nx_j$}\\
\scalebox{.9}{$=w_i(s_{-i})P(s)-w_i(s_{-i})\gamma^T\ltimes_{j=1}^{i-1}x_j\ltimes {\bf 1}_{k_i}\ltimes_{j=i+1}^nx_j$}\\
\scalebox{.9}{$=w_i(s_{-i})P(s)-w_i(s_{-i})\sum_{x\in S_i}P(x,s_{-i}).$}
\end{array}
$$
$2\Rightarrow 3:$ Plunging (\ref{e3.31}) into the left hand sides of (\ref{e3.32}) and (\ref{e3.33}) respectively, it is easy to verify these two equations.

$3\Rightarrow 1:$ According to Definition \ref{d1.1}, (\ref{e3.33}) shows that $G$ is a coset weighted potential game, then we only need to verify its orthogonality to ${\cal N}$ by using (\ref{e3.32}).
$$
\begin{array}{llc}
\scalebox{.9}{$\sum_{x\in S_i}c_i(x,s_{-i})=\sum_{x\in S_i}V^c_i\ltimes_{j=1}^{i-1}x_j\ltimes x\ltimes_{j=i+1}^nx_j$}\\
\scalebox{.9}{$=V^c_i\sum_{x\in S_i}\ltimes_{j=1}^{i-1}x_j\ltimes x\ltimes_{j=i+1}^nx_j$}\\
\scalebox{.9}{$=V^c_i\ltimes_{j=1}^{i-1}x_j\ltimes {\bf 1}_{k_i}\ltimes_{j=i+1}^nx_j$}\\
\scalebox{.9}{$=V^c_iE_iE_i^T\ltimes_{j=1}^nx_j=0.$}
\end{array}
$$
Then we have $\scalebox{.9}{$V^c_iE_iE_i^T=0$}$, which is equivalent to $\scalebox{.9}{$V^c_iE_i=0$}$, it follows that $\scalebox{.9}{$[V^c_1,V^c_2,\cdots,V^c_n]B^N=V_GB^N=0.$}$
Hence $G$ is is a coset weighted pure potential game, which is orthogonal to ${\cal N}$.
\hfill $\Box$
\begin{rem}\label{r3.2} We call (\ref{e3.61}) the geometric expression of coset weighted pure potential games, and (\ref{e3.31})-(\ref{e3.33}) are its algebraic expressions.
\end{rem}
\subsection{Subspace of coset weighted pure harmonic games ${\cal H}^{cw}$}
From the construction of $E_{cw}^{P^0}$, we have the dimension of coset weighted potential subspace as
$\scalebox{.9}{$\dim \left({\cal G}^{cw}_P\right)=k+\dsum_{j=1}^n\frac{k}{k_j}-1.$}$ Then the dimension of subspace ${\cal H}^{cw}$ is calculated as
\begin{align}\label{3.13}
\scalebox{.9}{$\dim \left({\cal H}^{cw}\right)=(n-1)k-\dsum_{j=1}^n\frac{k}{k_j}+1.$}
\end{align}
Set $\scalebox{.9}{$\psi^{cw}_n:=\left(\widetilde{E}_{cw}^P\right)^T$}$, obviously, we have
\begin{align}\label{3.13e}
\scalebox{.9}{${\cal H}^{cw}=\left(\widetilde{E}_{cw}^P\right)^{\perp}=\ker(\psi^{cw}_n).$}
\end{align}
Similar to the arguments in \cite{wang17}, we can construct
$$
\psi^{cw}_2=\begin{bmatrix}\begin{smallmatrix}
\Lambda_1&\Lambda_2\\
E_1^T\Lambda_1&0\\
0&E_2^T\Lambda_1
\end{smallmatrix}\end{bmatrix}.
$$
Set
$$
\begin{array}{cl}
x_{i_1,i_2}:=\begin{bmatrix}\begin{smallmatrix}
\Lambda_1^{-1}\left(\d_{k_1}^1-\d_{k_1}^{i_1}\right)\left(\d_{k_2}^1-\d_{k_2}^{i_2}\right)\\
~~\\
-\Lambda_2^{-1}\left(\d_{k_1}^{1}-\d_{k_1}^{i_1}\right)\left(\d_{k_2}^1-\d_{k_2}^{i_2}\right)\\
\end{smallmatrix}\end{bmatrix},\\
 \scalebox{.9}{$i_1=2,3,\cdots,k_1; ~i_2=2,3,\cdots,k_2.$}
\end{array}
$$
It is easy to see that
$$
x_{i_1,i_2}\in \ker(\psi^{cw}_2),\quad \scalebox{.9}{$i_1=2,3,\cdots,k_1; ~i_2=2,3,\cdots,k_2,$}
$$
and $\{x_{i_1,i_2}\;|\; i_1=2,3,\cdots,k_1; ~i_2=2,3,\cdots,k_2\}$ are linearly independent. From (\ref{3.13}), we calculate that $\dim({\cal H}^{cw}_2)=(k_1-1)(k_2-1)$, Hence, $\{x_{i_1,i_2}\;|\; i_1=2,3,\cdots,k_1; ~i_2=2,3,\cdots,k_2\}$ form a basis of ${\cal H}^{cw}_2$.

Next, we give an inductive method to construct $\psi^{cw}_n$.
\begin{lem}\label{l3.2} The matrix $\psi^{cw}_s$, $2\leq s\leq n$, can be recursively constructed by
\begin{align}\label{3.17}
	\psi^{cw}_p=
	\begin{bmatrix}\begin{smallmatrix}
	\psi^{cw}_{p-1} & \beta^{cw}_p\\
	{\bf 0}_{\frac{k}{k_p} \times(p-1)k} & (I_{\frac{k}{k_p}}\ot {\bf 1}_{k_p})\Lambda_1
	\end{smallmatrix}\end{bmatrix},
\end{align}
where $\beta^{cw}_p= [\Lambda_p,{\bf 0}_{k\times \frac{k}{k_1}},\cdots, {\bf 0}_{k\times \frac{k}{k_{p-1}}}]^T$, $k=\prod_{i=1}^{p}k_i$.
\end{lem}
According to Lemma \ref{l3.2}, it is easy to verify the following result by straightforward computations.
\begin{lem}\label{l3.3} If $x\in \ker(\psi^{cw}_{p-1})$, then
\begin{align}
\label{3.18}\begin{bmatrix}\begin{smallmatrix}
x \d_{k_p}^{i_{p}}\\
{\bf 0}_{k}
\end{smallmatrix}\end{bmatrix}\in \scalebox{.9}{$\ker(\psi^{cw}_p),\quad i_p=1,\cdots,k_p;$}\\
\notag
~\\
\label{3.19}\begin{array}{llc}
\begin{bmatrix}\begin{smallmatrix}
	\Lambda_1^{-1}(\d_{k_1}^1-\d_{k_1}^{i_1})\d_{k_2}^1\d_{k_3}^1 \cdots \d_{k_{p-1}}^1\\
	\Lambda_2^{-1}\d_{k_1}^{i_1}(\d_{k_2}^1-\d_{k_2}^{i_2})\d_{k_3}^1\cdots \d_{k_{p-1}}^1\\
	\vdots\\
	\Lambda_{p-1}^{-1}\d_{k_1}^{i_1}\d_{k_2}^{i_2}\d_{k_3}^{i_3}\cdots (\d_{k_{p-1}}^1-\d_{k_{p-1}}^{i_{p-1}})\\
    -\Lambda_p^{-1}\left(\d_{k_1}^1\d_{k_2}^1\cdots \d_{k_{p-1}}^1-\d_{k_1}^{i_1}\d_{k_2}^{i_2}\cdots \d_{k_{p-1}}^{i_{p-1}}\right)
	\end{smallmatrix}\end{bmatrix}\scalebox{.9}{$(\d_{k_p}^1-\d_{k_p}^{i_p})$}\\
\in \scalebox{.9}{$\ker(\psi^{cw}_p),~i_j=1,\cdots,k_j;~ j = 1, 2, \cdots, p.$}
\end{array}
\end{align}
\end{lem}
Define an index set as $\scalebox{.9}{$I=\{(i_1,\cdots,i_n)\;|\;i_p\in [1,~k_p]\}$}$, $p=1,\cdots,n$. Using Lemma \ref{l3.3}, we construct a set of vectors as follows, which are in $\ker(\psi_n)$.
$$
\begin{array}{llc}
\scalebox{.9}{$J_1$}=\left\{\left.\begin{bmatrix}\begin{smallmatrix}
			\Lambda_1^{-1}(\d_{k_1}^1-\d_{k_1}^{i_1})(\d_{k_2}^1-\d_{k_2}^{i_2})\d_{k_3}^{i_3} \cdots \d_{k_n}^{i_n}\\
			-\Lambda_2^{-1}(\d_{k_1}^1-\d_{k_1}^{i_1})(\d_{k_2}^1-\d_{k_2}^{i_2})\d_{k_3}^{i_3} \cdots \d_{k_n}^{i_n}\\
			{\bf 0}_{(n-2)k}
		\end{smallmatrix}\end{bmatrix}\;\right|\;\right.\\
~~~~~~~~~~~~~~~~~~~~\scalebox{.9}{$\left.i\in I, i_1\neq 1,i_2\neq 1\right\};$}
\end{array}
$$
$$
\begin{array}{llc}
\scalebox{.9}{$J_2$}=\left\{\left.\begin{bmatrix}\begin{smallmatrix}
\Lambda_1^{-1}(\d_{k_1}^1-\d_{k_1}^{i_1})\d_{k_2}^1(\d_{k_3}^1-\d_{k_3}^{i_3})\d_{k_4} ^{i_4}\cdots \d_{k_n}^{i_n}\\
\Lambda_2^{-1}\d_{k_1}^{i_1}(\d_{k_2}^1-\d_{k_2}^{i_2})(\d_{k_3}^1-\d_{k_3}^{i_3})\d_{k_4} ^{i_4} \cdots \d_{k_n}^{i_n}\\
-\Lambda_3^{-1}(\d_{k_1}^1\d_{k_2}^1-\d_{k_1}^{i_1}\d_{k_2}^{i_2})(\d_{k_3}^1-\d_{k_3}^{i_3})\d_{k_4} ^{i_4}\cdots \d_{k_n}^{i_n}\\
	{\bf 0}_{(n-3)k}
\end{smallmatrix}\end{bmatrix}\;\right|\;\right.\\
 ~~~~~~~~~~~~~~~~~~~~\scalebox{.9}{$\left.i\in I, (i_1,i_2)\neq {\bf 1}_2^T, i_3\neq 1\right\};$}
\end{array}
$$
$$
\vdots
$$
$$
\begin{array}{llc}
\scalebox{.9}{$J_{n-1}$}=\left\{\left.\begin{bmatrix}\begin{smallmatrix}
		\Lambda_1^{-1}(\d_{k_1}^1-\d_{k_1}^{i_1})\d_{k_2}^1\d_{k_3}^1\cdots \d_{k_{n-1}}^1 (\d_{k_n}^1-\d_{k_n}^{i_n})\\		
\Lambda_2^{-1}\d_{k_1}^{i_1}(\d_{k_2}^1-\d_{k_2}^{i_2})\d_{k_3}^1\cdots \d_{k_{n-1}}^1 (\d_{k_n}^1-\d_{k_n}^{i_n})\\
\vdots\\
\Lambda_{n-1}^{-1}		\d_{k_1}^{i_1}\d_{k_2}^{i_2}\d_{k_3}^{i_3}\cdots (\d_{k_{n-1}}^1-\d_{k_{n-1}}^{i_{n-1}}) (\d_{k_n}^1-\d_{k_n}^{i_n})	\\
-\Lambda_n^{-1}(\d_{k_1}^1\cdots \d_{k_{n-1}}^1-\d_{k_1}^{i_1}\cdots \d_{k_{n-1}}^{i_{n-1}} )(\d_{k_n}^1-\d_{k_n}^{i_n})			
\end{smallmatrix}\end{bmatrix}\;\right|\; \right.\\
 ~~~~~~~~~~~~~~~~~~~~\scalebox{.9}{$\left.i\in I, (i_1,\cdots,i_{n-1})\neq {\bf 1}_{n-1}^T, i_n\neq 1\right\}.$}
\end{array}
$$

Define a matrix as
\begin{align}\label{3.20}
\scalebox{.9}{$B_{cw}^H:=\left[J_1,J_2,\cdots,J_{n-1}\right].$}
\end{align}
Similar to the arguments in \cite{che16b,wang17}, we have the following result.
\begin{thm}\label{t3.3} $B_{cw}^H$ has full column rank and the subspace of coset weighted pure harmonic games is
\begin{align}\label{3.21}
\scalebox{.9}{${\cal H}^{cw}=\Span\left(B_{cw}^H\right),$}
\end{align}
where $\scalebox{.9}{$\Col(B_{cw}^H)$}$ is a basis of ${\cal H}^{cw}$.
\end{thm}
Next, we provide the algebraic expression of coset weighted pure harmonic games.
\begin{thm}\label{t3.6} Consider $G\in {\cal G}$. $G$ is a coset weighted pure harmonic game, if and only if, there exists a set of coset-depending weights $w_i(s_{-i})>0$, such that for any $s_{-i}\in S_{-i}$, and any $s_{-j}\in S_{-j}=\prod_{\substack{q\neq j\\q\neq 1 }}S_q$,
\begin{align}
\label{e3.34}\scalebox{.9}{$\sum^n_{i=1}w_i(s_{-i})c_i(s)=0;$}\\
\label{e3.35}\scalebox{.9}{$\sum_{\substack{x_j\in S_j\\  j\neq 1}}w_1(x_j,s_{-j})\sum_{x\in S_i}c_i(x,s_{-i})=0.$}
\end{align}
\end{thm}
\noindent{\bf Proof.} (Necessary) Assume the structure vector of $G$ is $\scalebox{.9}{$V_G=[V^c_1,V^c_2,\cdots,V^c_n].$}$ According to the orthogonality of (\ref{3.13e}), we have $\scalebox{.9}{$V_G^T\in \ker\left(\widetilde{E}_{cw}^P\right)^T.$}$ It follows that
\begin{align}\label{3.36}
\scalebox{.9}{$V_G\widetilde{E}_{cw}^P=\left[\sum^n_{i=1}V^c_i\Lambda_i,V^c_1\Lambda_1E_1,V^c_2\Lambda_1E_2,\cdots,V^c_n\Lambda_1E_n\right]=0.$}
\end{align}
Using Lemma \ref{l3.1}, Proposition \ref{p2.1} and \ref{p2.2}, we have
\begin{align}\label{e3.36}
\begin{array}{llc}
\scalebox{.9}{$\sum^n_{i=1}V^c_i\Lambda_i\ltimes_{j=1}^nx_j=\sum^n_{i=1}V^c_i(V^w_iE_i^T \otimes I_k)O_k^R\ltimes_{j=1}^nx_j$}\\
\scalebox{.9}{$=\sum^n_{i=1}V^w_iE_i^T(I_k\otimes V^c_i)\ltimes_{j=1}^nx_j\ltimes_{j=1}^nx_j$}\\
\scalebox{.9}{$=\sum^n_{i=1}V^w_iE_i^T\ltimes_{j=1}^nx_jV^c_i\ltimes_{j=1}^nx_j$}\\
\scalebox{.9}{$=\sum^n_{i=1}w_i(s_{-i})c_i(s)=0.$}
\end{array}
\end{align}
And we have
\begin{align}\label{e3.37}
\begin{array}{llc}
\scalebox{.9}{$V^c_i\Lambda_1E_iE_i^T\ltimes_{j=1}^nx_j=V^c_i\Lambda_1\ltimes_{j=1}^{i-1}x_j\ltimes {\bf 1}_{k_i}\ltimes_{j=i+1}^nx_j$}\\
\scalebox{.9}{$=V^c_i(V^w_1E_1^T \otimes I_k)O_k^R\ltimes_{j=1}^{i-1}x_j\ltimes {\bf 1}_{k_i}\ltimes_{j=i+1}^nx_j$}\\
\scalebox{.9}{$=V^w_1E_1^T\ltimes_{j=1}^{i-1}x_j\ltimes {\bf 1}_{k_i}\ltimes_{j=i+1}^nx_jV^c_i\ltimes_{j=1}^{i-1}x_j\ltimes {\bf 1}_{k_i}\ltimes_{j=i+1}^nx_j$}\\
\scalebox{.9}{$=\sum_{\substack{x_j\in S_j\\  j\neq 1}}w_1(x_j,s_{-j})\sum_{x\in S_i}c_i(x,s_{-i})=0.$}
\end{array}
\end{align}
(Sufficiency) It is clear that (\ref{e3.36}) and (\ref{e3.37}) can be deduced by (\ref{e3.34}) and (\ref{e3.35}) respectively. Because $\scalebox{.9}{$V^c_i\Lambda_1E_iE_i^T=0$}$ is equivalent to $\scalebox{.9}{$V^c_i\Lambda_1E_i=0$}$, hence, (\ref{e3.36}) and (\ref{e3.37}) can assure (\ref{3.36}), which leads to the conclusion.
\hfill $\Box$
\begin{rem}\label{r3.3} Theorem \ref{t3.3} gives the geometric expression of coset weighted pure harmonic games by providing $B_{cw}^H$. (\ref{e3.34}) and (\ref{e3.35}) are its corresponding algebraic expressions.
\end{rem}
\subsection{Numerical formulas of decomposed subspaces}
From the above arguments, the new orthogonal decomposition (\ref{1.5}) is established for every fixed coset-depending weight $w_i(s_{-i})$, $i=1,\cdots,n$. Then construct a basis matrix as $\scalebox{.9}{$B_{cw}:=[B_{cw}^{P^0},B^N,B_{cw}^H]$}$. Set
$\scalebox{.9}{$d_1=k-1$}$, $\scalebox{.9}{$d_2=\dsum_{j=1}^n k/k_j$}$, and $\scalebox{.9}{$d_3=(n-1)k-\dsum_{j=1}^n k/k_j+1$}$. Construct a set of coefficients as $\scalebox{.9}{$X^G_{cw}=[X_{cw}^P,~X^N,~X_{cw}^H]^T$}$, where $\scalebox{.9}{$X_{cw}^P\in \R^{d_1}$}$, $\scalebox{.9}{$X^N\in \R^{d_2}$}$ and $\scalebox{.9}{$X_{cw}^H\in \R^{d_3}$}$.
Assume the structure vector of $G$ is $V_G$, then we have
\begin{align}\label{4.3}
\scalebox{.9}{$(V_G)^T=B_{cw}X_{cw}^G=B_{cw}^{P^0}X_{cw}^P\oplus B^NX^N\oplus B_{cw}^HX_{cw}^H.$}
\end{align}
Using (\ref{4.3}), we can calculate all the decomposed components of a given game $G$ with fixed coset weights.
\begin{prp}\label{p4.1} Consider $G\in {\cal G}$. Let
\begin{align}\label{e6.1}
\scalebox{.9}{$\left[X_{cw}^P,~X^N,~X_{cw}^H\right]^T=(B_{cw})^{-1}(V_G)^T,$}
\end{align}
then
\begin{enumerate}
\item its coset weighted pure potential projection is:
\begin{align}\label{4.4}
\scalebox{.9}{$(V_G^{P^{cw}})^T=B_{cw}\left[X_{cw}^P~~0~~0\right]^T;$}
\end{align}
\item its coset weighted  pure harmonic projection is:
\begin{align}\label{4.5}
\scalebox{.9}{$(V_G^{H^{cw}})^T=B_{cw}\left[0~~0~~X_{cw}^H\right]^T;$}
\end{align}
\item its non-strategic projection is:
\begin{align}\label{4.6}
\scalebox{.9}{$(V_G^N)^T=B_{cw}\left[0~~X^N~~0\right]^T;$}
\end{align}
\item its coset weighted potential projection is:
\begin{align}\label{4.7}
\scalebox{.9}{$(V_{{\cal G}_P}^{cw})^T=(V_G^{P^{cw}})^T+(V_G^N)^T;$}
\end{align}
\item its coset weighted harmonic projection is:
\begin{align}\label{4.8}
\scalebox{.9}{$(V_{{\cal G}_H}^{cw})^T=(V_G^{H^{cw}})^T+(V_G^N)^T.$}
\end{align}
\end{enumerate}
\end{prp}
\begin{exa}\label{ex2} Recall Example \ref{ex1}, we can calculate all the decomposed components of $G$ with respect to $\scalebox{.9}{$V^w_1=[1,2]$}$ and $\scalebox{.9}{$V^w_2=[4,2]$}$. According to (\ref{2.2.3}), we have the structure vector of $G$ as $\scalebox{.9}{$V_G=[-1,2,0,3,3,3,5,4].$}$
Using (\ref{3.16}), it is easy to calculate that
$$
B_{cw}^P=\begin{bmatrix}\begin{smallmatrix}
\Lambda_1-\frac{1}{k_1}\Lambda_1E_1E_1^T\\
\Lambda_2-\frac{1}{k_2}\Lambda_2E_2E_2^T\\
\end{smallmatrix}\end{bmatrix}=\begin{bmatrix}\begin{smallmatrix}
 0.5  & 0  & -0.5  & 0\\
 0  & 1  & 0  & -1\\
 -0.5 & 0  & 0.5  & 0\\
 0  & -1  & 0  & 1\\
 2  & -2  & 0  & 0\\
 -2  & 2  & 0  & 0\\
 0  & 0  & 1  & -1\\
 0  & 0  & -1  & 1\\
\end{smallmatrix}\end{bmatrix}.
$$
 Construct $B_{cw}^{P^0}$ by deleting the last column of $B_{cw}^P$. According to (\ref{3.15}) and (\ref{3.20}), we have 
 $$
 \scalebox{.9}{$B^N$}=\begin{bmatrix}\begin{smallmatrix}
E_1&0\\
0&E_2\\
\end{smallmatrix}\end{bmatrix},
$$
and
$$
\begin{array}{ll}
\scalebox{.9}{$B_{cw}^H=[ 1  , -0.5  , -1  , 0.5  , -0.25  , 0.25 , 0.5 , -0.5]^T.$}
\end{array}
$$
Construct $\scalebox{.9}{$B_{cw}:=[B_{cw}^{P^0},B^N,B_{cw}^H].$}$
According to Proposition \ref{p4.1}, the coefficients can be calculated by (\ref{e6.1}) as
$\scalebox{.9}{$X_{cw}^P=[-0.5,~-0.5 ,~ 0.5]^T$}$, $\scalebox{.9}{$X^N=[-0.5,~2.5,~3,~4.5]^T$}$, and $\scalebox{.9}{$X_{cw}^H= 0$}$. Using formulas (\ref{4.4})-(\ref{4.8}), we have
$$
\begin{array}{lll}
\scalebox{.9}{$V_G^{P^{cw}}=[-0.5,~-0.5,~0.5,~0.5,~ 0,~0,~0.5,~-0.5];$}\\
\scalebox{.9}{$V_G^{H^{cw}}={\bf 0}^T_8;$}\\
\scalebox{.9}{$V_G^N=[-0.5,~2.5,~-0.5,~2.5,~3,~3,~4.5,~4.5];$}\\
\scalebox{.9}{$V_{{\cal G}_P}^{cw}=[-1,~2,~0,~3,~3,~3,~5,~4];$}\\
\scalebox{.9}{$V_{{\cal G}_H}^{cw}=[-0.5,~2.5,~-0.5,~2.5,~3,~3,~4.5,~4.5]=V_G^N.$}
\end{array}
$$
Hence, this game is a coset weighted potential game.
\end{exa}

\section{Conclusion}

This paper investigates a more general potential game, called a coset weighted potential game. Using the STP method, a coset weighted potential equation is presented to verify the coset weighted potential game, a corresponding formula is obtained to calculate the potential function. Some useful properties are explored. Finally, a new orthogonal decomposition of ${\cal G}$ with respect to fixed coset weights is obtained. The geometric and algebraic expressions of all the subspaces are provided respectively, and some formulas are given to calculate all the decomposed components.

\renewcommand\baselinestretch{0.5}

\end{document}